\documentclass[nopreprintline,1p,sort&compress]{elsarticle}
\usepackage[breaklinks=false,colorlinks=true,citecolor=black,urlcolor=blue,linkcolor=black,pdfpagemode=none,pdfstartview=FitV]{hyperref}
\pdfoutput=1
\usepackage{graphicx}
\usepackage{url}
\usepackage[dvipsnames]{xcolor}
\usepackage{algorithm} %
\usepackage[noend]{algpseudocode} %
\usepackage{amssymb,cite}
\usepackage{amsmath,amsfonts,mathrsfs}
\usepackage[inline]{enumitem}

\usepackage{bm} %
\usepackage{arydshln} %
\usepackage{cancel}

\allowdisplaybreaks

\newcommand{\bmat}[1]{\begingroup%
	\renewcommand*{\arraystretch}{1.1}
	\begin{bmatrix} #1 \end{bmatrix}
	\endgroup}

\newcommand{\blockcomment}[1]{}

\newcommand{\ceil}[1]{\left\lceil #1 \right\rceil}

\newcommand{\one}{\scalebox{1.05}{$\mathbf{1}$}}
\newcommand{\reals}{\ensuremath{\mathbb{R}}}

\newcommand{\integers}{\ensuremath{\mathbb{Z}}}

\newcommand{\naturals}{\ensuremath{\mathbb{N}}}
\newcommand{\stoch}{\mathcal{M}} %

\newcommand{\coloneq}{\mathrel{\mathop:}=}
\newcommand{\eqcolon}{=\mathrel{\mathop:}}
\newcommand{\given}{\mid}

\newcommand{\ddpart}[2]{\frac{\partial #1}{\partial #2}}

\newproof{proof}{Proof}

\newtheorem{theorem}{Theorem}
\newtheorem{Lemma}[theorem]{Lemma}

\newtheorem{Theorem}{Theorem}

\newtheorem{Proposition}[theorem]{Proposition}

\newcommand{\obs}{\mathcal{O}}

\begin{document}
\begin{frontmatter}
\title{Recovering Markov Models from Closed-Loop Data\tnoteref{t1}}
\tnotetext[t1]{Published in Automatica, May 2019, \href{https://doi.org/10.1016/j.automatica.2019.01.022}{10.1016/j.automatica.2019.01.022 }%
 --- \textcopyright 2019. This manuscript version is made available under the CC-BY-NC-ND 4.0 license (\href{http://creativecommons.org/licenses/by-nc-nd/4.0/}{creativecommons.org/licenses/by-nc-nd/4.0})}
\author[JESZ]{Jonathan Epperlein}\ead{jpepperlein(at)ie.ibm.com}\
\author[JESZ]{Sergiy Zhuk}\ead{sergiy.zhuk(at)ie.ibm.com}\
\author[RS]{Robert Shorten}\ead{robert.shorten(at)ucd.ie}
\address[RS]{University College Dublin, Dublin, Ireland}
\address[JESZ]{IBM Research, Dublin, Ireland}
\begin{abstract}
Situations in which recommender systems are used to augment
decision making are becoming prevalent in many application domains.
Almost always, these
prediction tools (recommenders) are created with a view to affecting behavioural
change. Clearly, successful
applications actuating behavioural change, affect the original model
underpinning the predictor, leading to an inconsistency. This feedback
loop is often not considered in standard machine learning
techniques which rely upon machine learning/statistical
learning machinery. The objective of this paper is to develop tools
that recover unbiased user models in the presence of recommenders. More specifically, we assume that we observe a time series
which is a trajectory of a Markov chain $\bm{R}$ modulated by
another Markov chain $\bm{S}$, i.e. the transition matrix of $\bm{R}$
is unknown and depends on the current state of
$\bm{S}$. The transition matrix of the latter is also
unknown. In other words, at each time instant, $\bm{S}$ selects a transition matrix for $\bm{R}$ within
a given set which consists of known and unknown matrices. The state of $\bm{S}$, in turn, depends on the
current state of $\bm{R}$ thus introducing a feedback loop. We
propose an Expectation-Maximization (EM) type algorithm,
which estimates the transition matrices of $\bm{S}$ and
$\bm{R}$. Experimental results are given to demonstrate the efficacy of the approach.
\end{abstract}
\end{frontmatter}

\section{Introduction}
\label{sect:intro}
Our starting point for this paper is a frequently encountered problem that arises in the Smart Cities domain. Many 
decision support/recommender systems that are designed to solve  Smart City problems are data-driven: that is 
data, sometimes in real time, is used to build models to drive the design of recommender systems. Almost 
always, these datasets are treated as if they were obtained in an open-loop setting, i.e.\ without recommender 
influence. However, this is rarely the case and frequently the effects of recommenders are inherent in datasets 
used for model building \citep{Lazer2014,Sinha_NIPS17,Shorten_IEEETech16,Bottou}. This creates new challenges for the 
design of decision support/recommender systems under feedback. In particular, as engineers, we must take into 
account the fact that when we make a prediction, then this prediction affects the behaviour of operators 
\citep{Cosley2003} and this, in turn, changes the data set upon which the original model was built.  Clearly, the 
aforementioned effect is related to classical closed-loop identification, which is itself a mature topic in both control 
and economics \citep{IV1,VanDenHof1995,Forssell1999}.

Notwithstanding this fact, and even though avoiding closed-loop effects in the design of recommender systems has been the subject of study in 
several Smart City applications \citep{Schlote2014b,Schlote2015}, the development of algorithms to identify {\em models}
in closed loop remains a challenging problem in the context of Smart Cities.  This is due to the fact that closed-loop questions that arise in Smart Cities are, for the most part, qualitatively different to those
arising in other areas. For example, in control theory closed-loop questions typically arise in the
context of deterministic and parametric models subject to noise, whereas in Smart Cities, typical problems are
 characterised by large-scale data sets that are generated by largely unknown stochastic processes.
 
 Our objective is to consider 
one such problem class that  arises in Smart City related research, where we seek to identify a user 
model, based on observations obtained when the user is acting under the influence of a recommender. 
A particular instance of such a problem arises in the automotive domain where drivers are characterised using Markovian 
 models \citep{Krumm2008,Epperlein2018}, but where observations are obtained under the influence of recommenders acting 
on the driver. Motivated by such applications, we seek to develop methods to account for the 
effect of these recommender systems in data sets.
More formally, we shall consider systems with the
following structure: the \emph{process}, which generates the data; the \emph{model},
which represents the behaviour of the process; and the \emph{decision support tool}, which intermittently influences the process. In our setup, data from the process is used to build the
model. Typically, the model is used to construct a decision support
tool which itself then influences the process directly. This creates a
feedback loop in which the process, decision support tool and the
model are interconnected in a complicated manner. As a result, the
effect of the decision support tool is to bias the data being
generated by the process, and consequently to bias any model that is
constructed naively from the data.

To provide a little more context, and to return to the automotive example, we now illustrate such effects by means of the following application
that we have developed in the context of our automotive
research\footnote{\url{https://www.youtube.com/watch?v=KUKxZZByIUM}}. Consider a driver who drives a car regularly. In order
to design a recommender system for this driver we would like to build
a model of his/her behaviour. For example, in order to warn the driver
of, say, roadworks, along a likely route, we might use this model to predict the route of
the driver. 
\begin{figure}[htb]
	\centering
		\includegraphics[width=0.4\columnwidth]{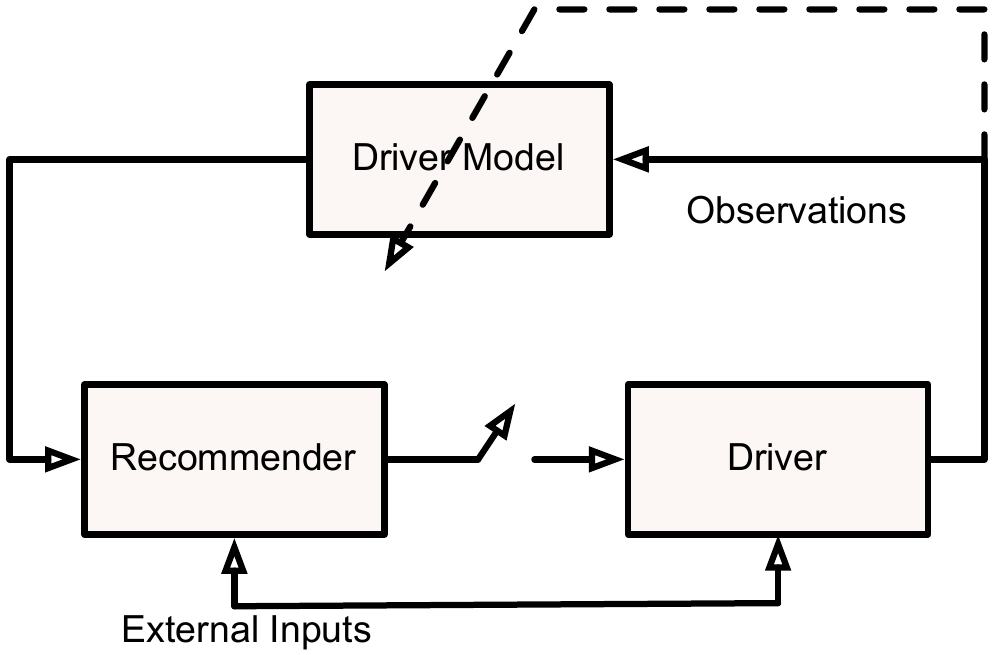}
		\caption{Automotive Recommender Architecture}
		\label{fig:arch}
\end{figure}
A schematic of the proposed in-car architecture is depicted in Figure \ref{fig:arch}. The recommender uses a model of driver behaviour to issue intermittent
recommendations. Observations of driver behaviour are then used to build a refined driver model which in turn is used as an input to the recommender system. Clearly,
the effect of the recommender is to bias the driver model over time,
thus eventually rendering the latter ineffective as an input to the recommender. The problems are exacerbated
in many practical systems due to the presence of several unknown
third-party recommender systems (Google Maps, Siri etc.), and  by the
fact that the {\em driver model} may operate from {\em {birth-to-death}}\protect\footnote{ by that we mean that the driver {\em always} operates under the potential influence of a recommender. Thus, given any observation, we do not know whether the recommender is acting, or the driver.}
in closed loop. This latter fact makes it difficult, or impossible, to even estimate an initial model of driver behaviour. Clearly, in such applications it
is absolutely necessary to develop techniques that
extract the behaviour of the driver while under the influence of the
feedback from a number of recommender systems. The results presented in this paper
represent our first small step in this direction.

\subsection{General Comments on Related Research Directions}
\label{sec:state-art}

Dealing with bias arising from closed-loop behaviour is a problem that has arisen in several application domains. In fact, in control theory, the related topic of
 closed-loop identification is considered to be a very mature area
 \citep{norton2009introduction,soderstrom1989system}. Roughly
 speaking, this topic is concerned with building models of dynamic
 systems while they are being regulated by a controller. A related
 scenario arises in some adaptive systems when the controller itself
 is being adjusted on the basis of the dynamical systems model. As in
 our example, the controller action will bias the estimation of the
 model parameters. While many established techniques in control theory
 exist for dealing with such effects, these typically exploit known
 properties of the process noise and an assumed model structure to
 un-bias the estimates. Typically, structures such as ARMAX models are
 assumed to capture the nature of the system dynamics. Recent work on
 {\em intermittent feedback} \citep{Gollee11} is also related in spirit
 to these approaches where control design techniques to deal with
 feedback loops that are intermittently broken are developed.  Before  proceeding it is worth noting that  {closed-loop effects have also been explored in the economics literature} \citep{HalVarian_NASUS17,MHE08, IV1,IV3}.
 
More recently, several authors in the context of Smart Cities and recommender systems~\citep{Shorten_IEEETech16,Sinha_NIPS17}, have realised that closed-loop effects represent a fundamental
challenge in the design of recommender systems. In \citep{Shorten_IEEETech16} the authors discuss the inherent
{\em closed-loop} nature of data-sets in cities, and in \citep{Sinha_NIPS17} 
 explicitly discuss the influence of feedback
on the fidelity of recommender systems.  {As an example of a specific
result, \citep{Sinha_NIPS17}  presents an empirical technique for collaborative filtering to recover
user rankings in the presence of a recommender under an assumed interaction model between user and
recommender, which is similar in spirit to the aforementioned problem
of constructing models from data possibly biased by a recommender, but
does not consider sequential models.}

The work reported here is closely
aligned with stochastic models, unlike most of the approaches outlined in the first paragraph
of Section~\ref{sec:state-art}. Specifically, in this paper we are interested in reconstructing Markov models that
are operating under the influence of a recommender. To this end we
assume that (i) \emph{recommenders and users can both be modelled by Markov
  chains}, and (ii) \emph{recommendations are either accepted fully or have no influence at all, i.e.\
  every decision is made by either the user or the recommender, never by a combination of the two}. We note here that the Markovian assumption of user and recommender behaviour
is convenient for many applications: for example, in the automotive domains~\citep{Krumm2008,JulienRoutePred}.

In this context, the present work is also related to classes of \emph{mixture} and \emph{latent variable models}, such as well-known hidden Markov (HMM) and mixture-of-experts (ME) models~\citep{JordanMarkovMixtures,JacobsMixtureExperts}. In the latter case, a latent Markov chain selects from a set of parametrisations of a visible process, however there is no closed-loop modulation (i.e.\  the modulated visible process is not allowed to in turn modulate the modulating latent process), and the visible process is static, subject to noise, not a Markov chain itself.
Our work is most related to~\citep{Ephraim_IEEETSP2009} but again with the important distinction of closed-loop modulation. A similar concept of regime switching time series
models is used in econometrics~\citep{Anderson2009-Lange}: these
models allow parameters of the conditional mean and variance to vary
according to some finite-valued stochastic process with states or
regimes. However, the observations are assumed to be
generated by a deterministic process with random noise, and the
latent (switching) process is either a Markov chain independent of the
past observations or is a deterministic function
of the past observations. In contrast, we introduce the closed-loop modulation as discussed above. Yet another related model is Markov jump linear systems, see e.g.~\citep{Costa2006}, where a (latent, autonomous) Markov chain selects the parameters of a (visible) dynamical system, whereas in our case, the visible part is a stochastic process on a discrete state space and it can modulate the latent process.

More specific technical comments to place our work in the context of reconstructing Markov
models from data, and a brief discussion of practical issues including \emph{identifiability} and convergence speed in terms of the number of samples are given below in Section~\ref{sec:HMM} after the formal description of the proposed closed-loop Markov-modulated Markov chain models.

\subsection{Preliminaries}\label{sec:prelim}
\textbf{Notation.}~To compactly represent discrete state spaces we write $[N]\coloneq\{1,\dotsc,N\}$. For a function $M:[N]\rightarrow\reals^{m\times n}$ mapping such a
discrete finite set to a set of matrices, we refer to each value
$M(k)\in\reals^{m\times n}$ as a \emph{page} of $M$. Matrices will be denoted by capital letters, their elements by the
same letter in lower case, and we denote the set of $n\times m$
row-stochastic matrices, i.e.\ matrices with non-negative entries such
that every row sums up to 1, by $\stoch^{n\times m}$, and
$\stoch^{m}\coloneq \stoch^{m\times m}$.
For compatible matrices, $M\otimes N$ is the Kronecker product and $M\circ N$ denotes the Hadamard (or element-wise) product.
A \emph{partition} $\Gamma$ of $[N]$ is a set $\{\Gamma_1,\dotsc,\Gamma_p\}$ such that $\Gamma_i\subseteq[N],\;\cup_i\Gamma=[N],\;\Gamma_i\cap\Gamma_j=\{\}$ $\forall i\neq j$. Each partition then also defines a 
membership function $\gamma: [N]\rightarrow [p]$ by $
	\gamma(i) \coloneq k \text{ such that } i\in\Gamma_k
$. 
We write $P(W=w)$ for the probability of the event that a
realisation of the discrete random variable $W$ equals $w$, and $P(W=w
\given V=v)$ the probability of that same event $W=w$ conditioned on
the event $V=v$. We shall
denote random variables by capital letters and, where appropriate, their realisations by
the same letter in lower case. For convenience we will sometimes write
$P(w \given v)$ instead of $P(W=w \given V=v)$ if there is no risk of
ambiguity, and, for a set of parameters $\mu$ parametrising a
probability distribution, $P(W=w \given \mu)$ is taken to denote the
probability of the event $W=w$ if the parameters are set to $\mu$.

\textbf{Markov chains} herein are sequences of random variables $\{X_t\}_t$ indexed by the time $t\in\naturals=\{0,1,2,\dotsc\}$. The realisation $x_t\in [N]$ of $X_t$ is the \emph{state} of the Markov chain at time $t$, and $[N]$ is its \emph{state space}. The probability distribution of $X_0$ is denoted by $\pi_0$, and the probability distribution of each following state is given by \(
 	P(X_t = j \given X_{t-1}=i) = a_{ij}
\); the matrix $A\in\stoch^{N}$ with entries $a_{ij}$ is the \emph{transition probability matrix}.

\section{Problem Statement  {and Model}}\label{sec:prob}
As noted above, we assume that the driver and the recommender are
Markovian. Given a possibly incomplete description of Markov chains
modelling the recommender systems, and no knowledge of when these
systems are engaged, our aim is to \emph{estimate the probability transition matrix of the Markov chain representing the
driver}, and the levels of engagement of each recommender, using only
observed data. In what follows we formalize this setup, and give an
expectation-maximization (EM) algorithm to estimate the parameters of the
unknown driver model. 
\subsection{``Open-Loop'' Markov-Modulated Markov Chains}\label{ssec:olmmmc}
Consider a Markov chain $\bm{R}$ with state space $[R]$ and state $r_t\in[R]$, in which the transition probabilities
\begin{equation}
	P(R_t=j\given R_{t-1}=i, S_{t-1}=s) = a^R_{ij}(s) \label{eq:which_AR}
\end{equation}
depend on a \textit{latent} random variable $S_t$. We can say that the
Markov chain is modulated by the random variable $S_t$, and if $S_t$
is itself the state of another Markov chain $\bm{S}$ with transition
matrix $A^S$ and state space $[S]$, then we are dealing with a
\textit{Markov-modulated Markov chain}; Markov modulation is an established model in the literature on  inhomogeneous stochastic processes, see e.g.~\citep{Ephraim_IEEETSP2009}.

Formally, the Markov-modulated Markov chain is defined by the tuple $\mu=(\pi^R, \pi^S,\allowbreak A^R(\cdot),A^S)$,
where $A^R: [S]\rightarrow \stoch^R$ and $\pi^S,\pi^R$ denote the distributions of $S_0$ and $R_0$, respectively. That means for instance that if $\bm{S}$ only has a single state $s$, $\mu$ is a regular Markov chain with transition matrix $A^R(s)$ and initial probability $\pi^R$. We assume that we observe the state of $\bm{R}$, but not the state of $\bm{S}$.
\\
Because the transition probabilities in the latent Markov chain $\bm{S}$ do not depend on the state of the visible chain $\bm{R}$, we refer to $\mu$ as an \emph{open-loop} Markov-modulated Markov chain (ol3MC) to distinguish it from what follows. This models the case when the switching between the transition matrices $A^{R}$ occurs independently of the current state $r_{t}$ of $\bm{R}$.

The joint process  ${Q_t}={(S_t, R_t)}$ has transition probabilities
\begin{multline*}
  P( Q_t = (s',r') \given Q_{t-1}=(s,  r) )
= \\
P(R_t=r' \given R_{t-1}=r, S_{t-1}=s, \cancel{S_t=s'})
\cdot P(S_t=s' \given \cancel{R_{t-1}=r,} S_{t-1}=s) 
=  a^R_{rr'}(s) \, a^S_{ss'},
\end{multline*}
where the first cancellation means that the decision at time $t-1$ is not influenced by the state of the modulating random variable at time $t$, and the second cancellation follows from the open-loop assumption, i.e. that the modulating Markov chain $\bm{S}$ evolves independently of $\bm{R}$. The estimation of $\bm{S}$ and $\bm{R}$ for the case of continuous time ol3MCs has been discussed in~\citep{Ephraim_IEEETSP2009}.

\textbf{Remark:} We are dealing with the case when the data consists of a finite time series
 of observations of a single trajectory  $(r_0r_1\dotsm r_T)$ of the Markov chain $\bm{R}$ and no (estimate of the) distribution of $R_t$ is available. While if the distributions are available, standard methods of state-space identification apply, here the estimation of the parameters of $\mu$ requires statistical methods such as maximum likelihood estimation.

\subsection{Closed-Loop Markov-Modulated Markov Chains}\label{ssec:clmmmc}
\begin{algorithm}[tb]
\begin{algorithmic}[1]
	\Procedure{cl3MC Trajectory}{$ \mu=(\pi^R,\allowbreak \pi^S,\allowbreak A^R, \allowbreak A^S;\allowbreak \Gamma)$, $T$ }
		\Statex \emph{Initialization}
			\State draw $s_0$ from $\pi^S$, $r_0$ from $\pi^R$ 
			\State	 $t \gets 0$
		\Statex \emph{Iteration}
			\While{$t\leq T-1$}
				\Statex \emph{Transition in $\bm{S}$}
					\State $\gamma' \gets \gamma(r_t)$ \Comment{active page in $\bm{S}$}
					\State draw $s_{t+1}$ from $\bmat{a^S_{s_{t}1}(\gamma')&\dotsm&a^S_{s_{t}S}(\gamma')}$
				\Statex \emph{Transition in $\bm{R}$}
					\State draw $r_{t+1}$ from $\bmat{a^R_{r_{t}1}(s_t)&\dotsm&a^R_{r_{t}R}(s_t)}$
			\EndWhile
			\State \textbf{return} $\obs=(r_0\dots r_T)$
	\EndProcedure
\end{algorithmic}
\caption{Trajectory generation for a given cl3MC.}
\label{alg:cl3MC}
\end{algorithm}
As a generalisation, we consider the case where $\bm{S}$ is dependent on the state of $\bm{R}$: that is, the probabilities of transitioning from one state $s$ to another state $s'$ then do depend on what the current state $r_t$ is. We will be referring to this as a \emph{closed-loop Markov-modulated Markov chain} or \emph{cl3MC} for short. A cl3MC can be used to model that one transition matrix $A^R(s)$ might be more likely to be switched to in some regions of the visible state space $[R]$, or that switching can only occur when the system is in specific configurations. This is exactly the situation which arises in our automotive example, see Section~\ref{ssec:driver}. 

Formally, we now also allow for the latent Markov chain $\bm{S}$ to be modulated by the current state of the visible chain $\bm{R}$. To keep the developments general, assume that -- instead of one page of $A^S$ corresponding to \emph{each }state of $\bm{R}$ -- there is a partition $\Gamma=\{\Gamma_1,\dotsc,\Gamma_p\}$ of $[R]$ such that there is a page in $A^S$ for each $\Gamma_i$. Hence, we now have $A^S:[p]\rightarrow\stoch^{S}$, with
\[
	P(S_t = j \given S_{t-1}=i, R_{t-1}=r) = a^S_{ij}( \gamma(r))\,.
\]
The open-loop case then corresponds to 
$\Gamma=\{[R]\}$ (i.e.\ $p=1$ and $\gamma(r)\equiv1$) and the joint process ${Q_t}={(S_t, R_t)}$  has transition probabilities (compare to the open-loop formula above):
\begin{multline}\label{eq:clMMMC}
  P( Q_t = (s',r') \given Q_{t-1}=(s,  r) )
= \\
 P(R_t=r' \given R_{t-1}=r, S_{t-1}=s ) \cdot P(S_t=s' \given R_{t-1}=r,\, S_{t-1}=s)\\
= 	 a^R_{rr'}(s) \, a^S_{ss'}(\gamma(r)).
\end{multline}
Such a cl3MC is represented by a tuple $\mu=(\pi^R,\allowbreak\pi^S,\allowbreak A^R(\cdot),\allowbreak A^S(\cdot);\Gamma)$, where now, $A^S(\cdot)$ has pages, too.

To further illustrate the operation of a cl3MC model, Algorithm~\ref{alg:cl3MC} details how a realization of the stochastic process described by it, i.e.\ a trajectory, is generated.

\subsection{{Relationship with Hidden Markov Models}}\label{sec:HMM}
There is a close relationship between closed-loop Markov modulated Markov chains and Hidden Markov Models (HMMs). Formally:
\begin{Proposition}\label{prop:HMM}
$\mu=(\pi^R,\pi^S,A^R(\cdot),A^S(\cdot);\Gamma)$ defines the same visible process $\{R_t\}$ as the Hidden Markov Model $\lambda=(\pi, W, B)$, with
\begin{equation*}
	\begin{split}
		\pi & = \pi^S\otimes\pi^R 	, \qquad
		B = 1_S\otimes I_R , \\[1ex]
		w_{ij} &= a^S_{\ceil{i/R}\,\ceil{j/R}}\bigl(\gamma(j_R(i))\bigr) \, a^R_{j_R(i) \, j_R(j)} \bigl( \ceil{i/R} \bigr),
	\end{split}
\end{equation*}
    	where $i,j=1,\dotsc,RS$, hence $W\in\stoch^{RS,RS}$.
\end{Proposition}
 {Here,  $j_R(k)\coloneq (k-1) (\bmod R) + 1$ and $\ceil{p} \coloneq \inf \{ k \in\integers \,\mid\, k\geq p \}$.
The proof of Proposition~\ref{prop:HMM} is mainly algebraic and can be found in~\citep{clMMMC_arxiv}.}

\textbf{Remark:} While Proposition~\ref{prop:HMM} maps a given cl3MC to an HMM which from the outside looks the same as, this mapping is not reversible: not every HMM represents a cl3MC, and most importantly, parameter estimation algorithms such as the standard Baum-Welch algorithm can not be used to estimate the parameters of a cl3MC, because they do not ``respect the structure'' of the matrix $W$: the HMM $\lambda$ is defined by $(RS)^2-RS+(RS-1)$ free parameters (the entries of $W$ and the entries of $\pi$ with the stochasticity constraints taken into account), whereas the corresponding cl3MC requires only $(pS-1)(S-1)+(SR-1)(R-1)$ parameters\footnote{
	Note that $(pS-1)(S-1)+(SR-1)(R-1) < (RS)^2-1$ %
	for $R+S>2$.}.
Hence, it is not possible to estimate the parameters of $\lambda$ and then compute the ones of $\mu$; instead, we develop an EM-algorithm to estimate the parameters of $\mu$ directly in Section~\ref{ssec:EM}.

\textbf{Identifiability: } \label{par:identifiability} Given the close relationship between {cl3MCs}
and HMMs outlined above, one should expect that identifiability issues
for cl3MCs bear close resemblance to those of HMMs. By
identifiability we mean the following: assume that $(r_1\dots r_T)$ has been generated by the ``true model''
$\mu^{\mathrm{true}}$; under which conditions and in what sense will the
estimate $\mu^{\mathrm{est}}$
 converge to $\mu^{\mathrm{true}}$ if $T\to\infty$? For HMMs this question was partially answered in~\citep{PetriAMS1969}, namely it was shown that there is an open, full-measure subset $U$ of all HMMs, such that the sequence of estimates of the BW algorithm converges to $\lambda^{\mathrm{true}}$ (or a trivial permutation of it), provided the starting model is chosen within $U$, $\lambda^{\mathrm{true}}_i>\delta>0$ and $T\to\infty$. However, the structure of $U$ and convergence speed in terms of the number of samples were not described, and, to the best of our knowledge, these questions are still open.

For cl3MCs, similarly and trivially, any permutation of $A^S$ and the corresponding pages of $A^R$, which amounts to relabelling the hidden states $s$, yields the same visible process.
 However, there are examples of sets of HMMs $\lambda$, which are not permutations of each other, yet generate the same observable process; see~\citep{Blackwell1957,GilbertAMS1959}. Interestingly, those examples involve the special case of \emph{partially observable} Markov chains, a subclass of HMMs with emissions matrices $B$ having entries that are either $1$ or $0$. Comparing to Proposition~\ref{prop:HMM}, a cl3MC has close correspondence to an HMM of this class.
 This suggests that, in practice, the set of maximisers of the likelihood may be wider than the aforementioned set of permutations of $\mu^{\mathrm{true}}$; our numerical experiments in Section~\ref{sec:exm} also suggest that, in general, we cannot recover the true model $\mu^{\mathrm{true}}$, even up to trivial permutations, from observing only trajectories of $\bm{R}$.
 However, in the case of partial knowledge of elements of $A^R$ we can recover $A^S$ and the unknown portion of $A^R$. 
 Hence, for the  ``driver-recommender'' problem the proposed method is of practical value. Estimates of the minimum amount of prior knowledge necessary are the subject of future research.
\section{Likelihood and Parameter Estimation}\label{sec:algo}

In this section we develop an iterative algorithm to estimate the parameters of a cl3MC $\mu=(\pi^R,\allowbreak \pi^S,\allowbreak A^S(\cdot),\allowbreak A^R(\cdot); \Gamma)$ given a sequence of observations $(r_0r_1r_2\dotsm r_T)$, a partition $\Gamma=\{\Gamma_1,\dotsc,\Gamma_p\}$ of $[R]$ and the size $S$ of the state space of $\bm{S}$.
The derivation is close in spirit to the classical Baum-Welch (BW)
algorithm (see e.g.~\citep{Rabiner1989} and the numerous references
therein): our algorithm maximises at every iteration a lower bound on
the likelihood improvement, and gives rise to re-estimation
formulae~\eqref{eq:reest} that utilise forward and backward variables
which differ in subtle ways from the ones of the BW algorithm. 

\subsection{Likelihood of $\mu$, Forward- and Backward Variables}
\label{ssec:likely}
Since the estimate to be obtained is a maximum likelihood (ML)
estimate, the efficient computation of the likelihood of a given
cl3MC $\mu$ plays a central role in what follows. For a given $\mu$, the joint probability of
sequences $(r_0r_1\dotsm r_T)$ and $(s_0s_1\dots s_T)$
being the trajectories of the visible Markov chain $\bm{R}$ and latent Markov chain $\bm{S}$ is
\begin{equation}
  \label{eq:Prs}
  \begin{split}
  P(s_0,&\dotsc,s_T,r_0,\dotsc,r_T\given \mu)\\
  &=\pi^R_{r_0}\pi^S_{s_0}\prod_{t=1}^T P(s_t,r_t\given
  s_{t-1},r_{t-1},\mu)
  = \pi^R_{r_0}\pi^S_{s_0}\prod_{t=1}^T a^S_{s_{t-1} s_t} (\gamma(r_{t-1})) a^R_{r_{t-1} r_t} (s_{t-1})
\end{split}
\end{equation}
where the last equality follows by~\eqref{eq:clMMMC}. This allows us
to compute the probability of observing a sequence $(r_0r_1\dotsm r_T)$ given $\mu$ as follows:
\begin{equation}
\label{eq:lmu}
  \begin{split}
P&(r_0r_1\dotsm r_T\given \mu)\\
&= \sum_{s_0\in[S] \dotsb s_T\in[S]}P(s_0,\dotsc,s_T,r_0,\dotsc,r_T\given \mu)\\
&= \pi^R_{r_0}\sum_{s_0\dotsb s_T}\pi^S_{s_0}\prod_{t=1}^T a^S_{s_{t-1} s_t} (\gamma(r_{t-1})) a^R_{r_{t-1} r_t} (s_{t-1})
\eqcolon\ell(\mu)
              \end{split}
            \end{equation}
where $\mu\mapsto \ell(\mu)$ is the likelihood of the model $\mu$.
Computation using this direct expression requires on the order of $2\times T\times S^T$ operations, and is hence not feasible for large $T$. Instead, we define the \textit{forward variable} $\alpha_t$ with $S$ elements
\begin{equation}
	\alpha_t(i) \coloneq  P(S_t=i, R_0=r_0, \dotsc, R_t=r_t\given \mu)
\end{equation}
which can be computed iteratively as follows: $\alpha_0(j) = \pi^S_j \pi^R_{r_0}$ and 
\[
\alpha_t(j) = \sum_{i=1}^S \alpha_{t-1}(i) a^S_{ij}(\gamma(r_{t-1})) a^R_{r_{t-1}r_t}(i)
, \quad j=1,\dotsc S,
\] or, in matrix form: $\alpha_0 = \pi^S \pi^R_{r_0}$ and
\begin{equation}
  \label{eq:alpha}
  \alpha_t = \left(A^S(\gamma(r_{t-1}))\right)^\intercal\,\left( a^R_{r_{t-1}r_t}(:) \circ \alpha_{t-1}  \right),
\end{equation}
where the notation $a^R_{r_{t-1}r_t}(:)$ means a column vector of the $(r_{t-1},r_t)$-elements of the matrix $A^R(k)$ as $k$ runs from $1$ to $S$.\footnote{Very much analogous to Matlab's \texttt{colon} notation, or slicing in \texttt{numpy}.}
It follows that
\begin{equation*}
	\ell(\mu) \!=\! \sum_{i=1}^S \!P(S_T=i, r_0,\dotsc,r_T \given \mu) \!= \!\sum_i \alpha_T(i) = \one^\intercal \alpha_T
\end{equation*}
can be computed with on the order of $TS^2$ computations.

An analogous concept that will be required later is the \textit{backward variable}
\begin{equation}
	\beta_t(i) \coloneq P(r_{t+1},\dotsc,r_T \given S_t=i, r_t),
\end{equation}
which can also be computed via iteration: $\beta_T(j) = 1 $, 
\[\beta_{t-1}(j) = \sum_i \beta_t(i) a^R_{r_{t-1}r_t}(j) a_{ji}(\gamma(r_{t-1})),
\]
or in matrix form: $\beta_T  = \one$ and
\begin{equation}\label{eq:beta}
  \beta_{t-1}= \left( A^S(\gamma(r_{t-1})) \beta_t \right) \circ a^R_{r_{t-1}r_t}(:).
\end{equation}
\subsection{Auxiliary Function $Q(\mu,\mu')$}
\label{ssec:Q}
Let $\Lambda$ denote the set of all cl3MCs. $\Lambda$ is then bounded and convex if we define convex combinations of cl3MCs $\mu=(\pi^R,\allowbreak\pi^S,\allowbreak A^R(\cdot), \allowbreak A^S(\cdot);\Gamma)$ and $\nu=(\rho^R,\allowbreak\rho^S,\allowbreak B^R(\cdot), \allowbreak B^S(\cdot);\Gamma)$ as
\begin{multline*}
\alpha\mu+(1-\alpha)\nu =\Bigl(\alpha \pi^R+(1-\alpha)\rho^R,\alpha \pi^S+(1-\alpha)\rho^S,\\
\alpha A^R+(1-\alpha)B^R,\alpha A^S+(1-\alpha)B^S;\Gamma\Bigr).
\end{multline*}
See~\citep{clMMMC_arxiv} for details.
Following~\citep{Baum_AnnalsMatStat1970}, we define the auxiliary function $Q(\mu,\mu')$ of $\mu,\mu'\in\Lambda$ by
\begin{equation}
	\label{eq:Q}
		Q(\mu,\mu') \coloneq \sum_{s_0\dotsb s_T} P  (s_0,\dotsc,s_T,r_0,\dotsc,r_T\given \mu) \cdot 
			\log P(s_0,\dotsc,s_T,r_0,\dotsc,r_T\given \mu'),
\end{equation}
where $s_0,\dotsc,s_T$ run through the $S^T$ possible sequences of the latent state $S_t$.
If parameters $\mu'_i$ are zero where $\mu_i>0$, then we can have the case $P(s_0,\dotsc,r_T\given \mu')=0$ and $P(s_0,\dotsc,r_T\given \mu)>0$; in this case $Q(\mu,\mu')\coloneq-\infty$. If $P(s_0,\dotsc,r_T\given \mu)=P(s_0,\dotsc,r_T\given \mu')=0$, we set $Q(\mu,\mu')=0$ which amounts to setting $0\log(0)=0$.\\
The following lemma establishes a representation for $Q$ in terms of the elements of $\mu'=(\pi^{R'},\pi^{S'},A^{S'},A^{R'};\Gamma)$:
\begin{Lemma}\label{lem:Q}
 	The function $Q(\mu,\mu')$ can be rewritten as
        \begin{multline}
          \label{eq:Qsum}
            	Q(\mu,\mu') =  \log\pi^{R'}_{r_0}\ell(\mu) +
                \sum_{i=1}^S \log \pi^{S'}_i \sum_{j=1}^S \xi_1(i,j) 
                +{\sum_{i,j=1}^S \sum_{t=1}^T  L_{ij}(r_{t-1},r_t) \xi_t(i,j)}
        \end{multline}
	where {$L_{ij}(m,n)\coloneq \log a^{S'}_{ij}(\gamma({m})) + \log a^{R'}_{mn} (i)$} and
	\begin{equation}\label{eq:xidef}
		\xi_t(i,j)\coloneq P(S_{t-1}=i,S_t=j,r_0,\dotsc,r_T \given \mu)
	\end{equation}
	can be computed as follows: 
	\begin{equation}\label{eq:xi}
		\xi_t(i,j) = \alpha_{t-1}(i) a^R_{r_{t-1}r_t}(i) a^S_{ij}(\gamma(r_{t-1})) \beta_t(j),
	\end{equation}
	and the variables carrying a $\bullet'$ constitute $\mu'$.
 \end{Lemma}
The proof is relegated to~\ref{sec:proofs}. Additionally, an application of Jensen's inequality yields (see~\citep{clMMMC_arxiv}):
\begin{Lemma}\label{thm:Q}
	The improvement in log-likelihood satisfies the lower bound
        \begin{equation}
          \label{eq:Ql}
	\ell(\mu) \bigl(\log\ell(\mu') - \log\ell(\mu)\bigr) \geq {Q(\mu,\mu')-Q(\mu,\mu)}\,.
        \end{equation}
      \end{Lemma}
\subsection{EM-Algorithm for Parameter Estimation}\label{ssec:EM}
 The algorithm proceeds by maximising the lower bound on the log-likelihood improvement set forth in~\eqref{eq:Ql} at every iteration. 
 
It should be clear from~\eqref{eq:Qsum} that the best estimate of
$\pi^R$ is the $r_0$-th canonical Euclidian basis vector $e_{r_0}$. The remaining parameters of $\mu$ can be iteratively estimated by repeatedly applying the following theorem:
\begin{Theorem}\label{thm:reest}
	The unique maximizer $\mu'=M(\mu)$ of $Q(\mu,\cdot)$ is given by
	\begin{subequations}\label{eq:reest}
	\begin{align}
		\pi^{S'}_i &= \frac{\sum_{j=1}^S \xi_1(i,j)}{\one^\intercal \alpha_T} 	\label{eq:reestpi}\\ %
		a^{S'}_{ij}(l) &= \frac{\sum_{t: \gamma(r_{t-1})=l} \xi_t(i,j)}{\sum_{k=1}^S \sum_{t: \gamma(r_{t-1})=l} \xi_t(i,k)} 	\label{eq:reestas}\\
		a^{R'}_{mn}(i) &= \frac{\sum_{k=1}^S\sum_{t:
                                 r_{t-1}=m, r_t=n}\xi_t(i,k)}  {\sum_{\nu=1}^R\sum_{k=1}^S
                                 \sum_{t: r_{t-1}=m, r_t=\nu} \xi_t(i,k)}, 		\label{eq:reestar}
	\end{align}
	\end{subequations}
	where $i,j=1,\dotsc,S$, $l=1,\dotsc,p$, and $m,n=1,\dotsc,R$.
\end{Theorem}
The proof is given in~\ref{sec:proofs}. Formulae~\eqref{eq:reest}
  provide the basis for the EM-type parameter estimation algorithm
  for cl3MC $\mu$: in its $k$-th iteration, the E-step consists of
  computing $\xi_t(\cdot,\cdot)$ from the current estimate $\mu^k$,
  and the M-step yields an updated estimate $\mu^{k+1}=M(\mu^k)$ with
  improved likelihood, see also the pseudocode in Algorithm~\ref{alg:EM} in the Appendix. Note that $M$ has a unique fixed point,
  $\mu^\infty=\lim_k\mu^k$ which is, at the same time, a stationary point (possibly a local maxima) of the likelihood,
  see~\citep{clMMMC_arxiv} for details.
 \section{Examples}\label{sec:exm}
Here we illustrate the algorithm's efficacy in two scenarios: first with synthetic data, i.e.\ data generated from a cl3MC, denoted $\mu^\text{true}$; and second, in a toy example of a practical application, estimation of driver behaviour. In both cases, we assume that one decision-maker, specifically the matrix $A^R(2)$, and $\Gamma$ are known.
 For implementation details, in particular how to avoid arithmetic underflow by scaling, and a pseudocode, see~\ref{ssec:scaling} and Algorithm~\ref{alg:EM} therein, and~\citep{clMMMC_arxiv}; for details on the experimental procedures, see~\ref{app:exp}.
\subsection{Synthetic Data}\label{ssec:synth}
To explore the relationship between estimation error and number of samples, repeated the following for several values $T\in [500,75000]$:
$N_e=100$ cl3MCs with $R=20$, $S=2$ and $\Gamma=\{[R]\}$ (i.e.\ the open-loop case) were generated, and then a trajectory of length $T$ for each of them. We then ran the algorithm with random initial guesses  $\pi^S$, $A^S$, and $A^R(1)$. The same was repeated for the same $N_e$ cl3MCs, only that now, $\Gamma$ was a randomly selected partition of order 2, so that $p=2$, and a second random page $A^S(2)$ was added to $A^S$. In both cases, we assume $A^R(2)$ and $\Gamma$ to be known. The modification to the algorithm is trivial: $A^R(2)$ is simply not re-estimated.

\begin{figure}[tbh]
\centering\includegraphics[keepaspectratio,width=.7\columnwidth]{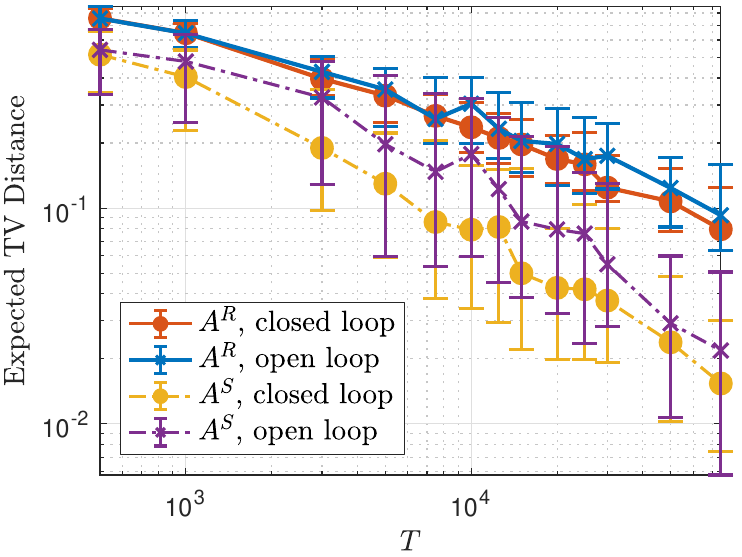}
\caption{Length of training sequen\-ces $T$ vs $\| A^{\bullet,\text{est}}(1)-A^{\bullet,\text{true}}(1) \|_\text{exp}$, defined in~\eqref{eq:expTV}. Shown are the medians and whiskers for the quartiles.}
\label{fig:err_vs_T}
\end{figure}

\begin{figure}[tbh]
\centering\includegraphics[keepaspectratio,width=.8\columnwidth]{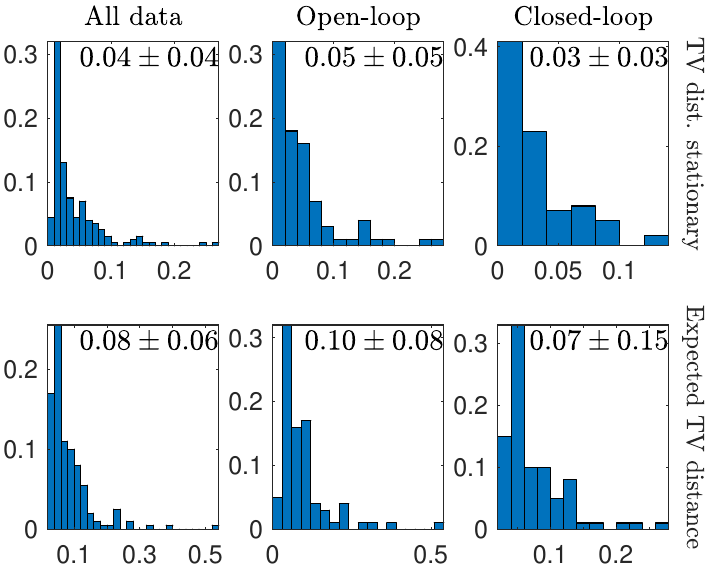}
\caption{The upper row shows the distribution of $\| A^{R,\text{est}}(1)-A^{R,\text{true}}(1) \|_\text{stat}$, see~\eqref{eq:statTV}, first for all $2 N_e=200$ pairs of $\mu^\text{true}$ and $\mu^\text{est}$, and then for the $N_e$ open-loop and $N_e$ closed-loop cases. In the bottom row, the same is shown for the metric $\| A^{R,\text{est}}(1)-A^{R,\text{true}}(1) \|_\text{exp}$, defined in~\eqref{eq:expTV}.}
\label{fig:page_fixed_bars_AR}
\end{figure}

As illustrated in Figures~\ref{fig:err_vs_T} and \ref{fig:page_fixed_bars_AR}, we recover $A^S$ and $A^R(1)$ to high accuracy for large enough $T$. ``Accuracy'' is hereby measured through statistical distances: since the transition matrices of Markov chains consist of probability distributions -- row $i$ being the distribution of the state following $i$ -- absolute or relative matrix norms are not a good measure of distance between Markov chains. Instead, we consider a statistical distance between the estimated and true probability distributions. One of the simplest such distances is the total variation (TV) distance (see e.g.~\citep[Ch.~4]{Levin2009}), which is given by the maximal difference in probability for any event between two distributions. For probability distributions $f$ and $g$ over a discrete space $\Omega$, this is simply
\begin{equation*}
	\|f - g\|_{TV} = \max_{A\subseteq\Omega} f(A) - g(A) = \frac{1}{2} \sum_{\omega\in \Omega} |f(\omega) - g(\omega)|.
\end{equation*}
We consider here two applications of TV distance to Markov chains. The first is to take the TV distance between the stationary distributions, which concretely amounts to considering the subset $\rho$ of the state space $[R]$ such that $P(X_t\in\rho \given A^{R,\text{est}}(1) ) - P(X_t\in\rho \given A^{R,\text{true}}(1) )$ is maximised (for large enough times $t$ such that the stationary distribution is reached). If we let $ \psi^\text{true}$ and $\psi^\text{est}$ denote the stationary distributions, then
\begin{equation}\label{eq:statTV}
	\| A^{R,\text{est}}(1)-A^{R,\text{true}}(1) \|_\text{stat} \coloneq \| \psi^\text{est} - \psi^\text{true} \|_{TV}.
\end{equation}
However, this is a coarse measure: different Markov chains can have equal stationary distributions. Hence, the second metric incorporates the distance between the individual rows by considering the expectation (under the true stationary distribution $\psi^\text{true}$) of the TV distance between the estimated and the true row; this equals the sum of the distances between the true and estimated transition probabilities from all states $i$, weighted by the probability of being in state $i$:
\begin{multline}\label{eq:expTV}
	\| A^{R,\text{est}}(1)-A^{R,\text{true}}(1) \|_\text{exp} \coloneq  \sum_i \psi^\text{true}_i \|A^{R,\text{est}}_{i:}(1)-A^{R,\text{true}}_{i:}(1)\|_{TV},
      \end{multline}
where $M_{i:}$ denotes the $i$-th row of matrix $M$. 

The effect of $T$ on the accuracy is explored in Figure~\ref{fig:err_vs_T}. The error appears to decay as a power of $T$, however this is simply an observation; a theoretical analysis of the sample complexity and decay rates is part of future work to be done.

For a representative value of $T=5\cdot10^{4}$, Figure~\ref{fig:page_fixed_bars_AR} drills down further into the experimental results; the distance for both introduced metrics is often below 10\%, but we also observe severe outliers. Note that we show $A^R(1)$ only, the analysis and results for $A^S(\cdot)$ are analogous and are hence omitted.

\subsection{A Model of Driver Behaviour}\label{ssec:driver}
Recent research, e.g.~\citep{Epperlein2018,JulienRoutePred,SimmonsBrowningZhangEtAl2006,Krumm2008}, suggests that Markov-based models are good approximations of driver behaviour and can be used e.g.\ for route prediction. Here, we illustrate how cl3MCs can be used to identify a driver's preferences when some trips are planned by a recommender system, whose preferences are known, while the other trips are planned by the driver.

\begin{figure}[tbh]
\parbox[c][\height][c]{.7\columnwidth}{
\includegraphics[keepaspectratio,width=.68\columnwidth]{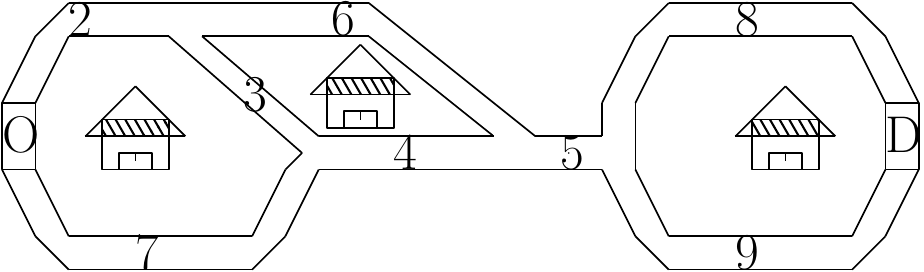}}\hfill
\parbox[c][\height][c]{.29\columnwidth}{
\includegraphics[keepaspectratio,width=.28\columnwidth]{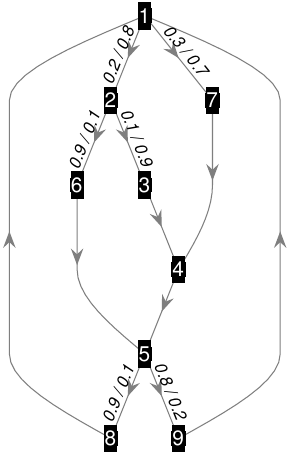}}
\caption{The map of our small toy model and its abstraction as the line graph of the road model. The direction of traffic is from left to right only. The origin and destination are merged into node 1. The weights denote the transition probabilities for the driver and the recommender system. When there is no weight, then the transition probability is 1 for both.}
\label{fig:schoolmap}
\end{figure}

Specifically, consider the map in the left panel of Figure~\ref{fig:schoolmap}, which depicts a (very small toy) model of a driver's possible routes from origin ``O'' to destination ``D.'' The houses, as an example, correspond to schools, that should be avoided in the hour before classes start and after classes end for the day, so there is a route past them and one around them. We assume that if a trip falls into that time frame, the recommender takes over and, with known probabilities, routes the driver either past or around each school; these probabilities make up $A^R(2)$. Otherwise, the driver follows his/her preferences, which constitute $A^R(1)$; this is the matrix we would like to estimate.
We generated $N_e=50$ sets of data by simulating $N_t=200$ trips on the graph shown in the right panel of Figure~\ref{fig:schoolmap}; this is the line graph of the map, where each road segment corresponds to a node, and an edge goes from node $i$ to node $j$ iff it is possible to turn into road segment $j$ from $i$. Each trip has a probability of $p_r=0.3$ to be planned by the recommender. If a trip was planned by the recommender (resp.\ driver), a trajectory was generated by a Markov chain with transition matrix $A^R(2)$ (resp.\ $A^R(1)$) originating in node 1 and terminating when returning to node 1.

\begin{figure}[tbh]
\centering
\includegraphics[keepaspectratio,width=.5\columnwidth]{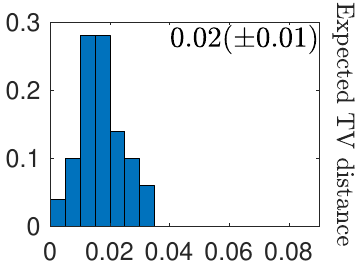}
\caption{Results from estimating the driver preferences in Section~\ref{ssec:driver}. Assuming open loop yields acceptable results; adding the information that only full trips are planned by either driver or recommender improves accuracy considerably. Mean and standard deviation are given.}
\label{fig:driver_results}
\end{figure}

For estimation in the cl3MC framework, all trips are then concatenated to form one long trajectory and $A^S$ and $A^R(1)$ are estimated for an ol3MC, i.e.\ for $\Gamma=[9]$. $A^{R,\text{est}}(1)$ is then an estimate of the driver preferences. The results are shown in  the first column of Figure~\ref{fig:driver_results} and are satisfactory already; however, we can leverage the \textit{closed}-loop framework to include the additional knowledge that the the decision maker (i.e.\ the page of $A^R$ used) can only change after a trip is finished. Because the decision which page of $A^R$ to use at time $t$ is made at $t-1$, see~\eqref{eq:which_AR}, this means we have to allow for the state of $\mathbf{S}$ to change on the road segments \emph{prior} to reaching the destination. We hence let $\Gamma = \bigl\{ \{8,9\}, \{1,\dotsc,7\} \bigr\}$ and $A^S(2)=I_2$. $A^S(1)$ needs to be identified. The results are shown in Figure~\ref{fig:driver_results}.
Additionally, we can interpret $\psi_2^{S,\text{est}}$, the second element of the stationary distribution of $A^{S,\text{est}}(1)$ as an estimate of $p_r$. For the open-loop case, we obtain $0.352(\pm0.0583)$, whereas the cl3MC estimation yields $0.295(\pm0.044)$.

\section{Concluding Remarks}
We consider the identification of user models acting under the influence of one or more recommender systems. %
As we have already discussed, actuating behavioural change affects the original model underpinning the predictor, leading
to an biased user models. Given this background, the specific contribution of this paper is to develop %
techniques in which unbiased estimates of user behaviour can be recovered in the case where recommenders, users, and switching between them can be parameterised in a Markovian manner, and where users and recommenders form part of a feedback system. Examples are given to present the efficacy of our approach.

{\small \textbf{Acknowledgements: }The authors would like to thank
Ming-Ming Liu and Yingqi Gu (University College Dublin) for help with
the numerical experiments, and Giovanni Russo and Jakub Mare\v{c}ek
(IBM Research Ireland) for valuable discussions.}

{\small
This work has been conducted within the ENABLE-S3 project that has
received funding from the ECSEL joint undertaking under grant
agreement NO 692455. This joint undertaking receives support from the
European Union's HORIZON 2020 Research and Innovation programme and
Austria, Denmark, Germany, Finland, Czech Republic, Italy, Spain,
Portugal, Poland, Ireland, Belgium, France, Netherlands, United
Kingdom, Slovakia, Norway. \\
Robert Shorten was also partially supported 
by SFI grant 16/IA/4610.}

\appendix
\section*{Appendix}
\section{{Scaling Issues}}\label{ssec:scaling}
{
Since the computations of $\alpha_t$ and $\beta_{T-t}$ according to~\eqref{eq:alpha} and \eqref{eq:beta} involve
multiplications of on the order of $2t$ numbers less than 1, for large
$T$, they will be close to, or below, machine precision. The
re-estimation~\eqref{eq:reest} then requires division of very small
numbers, which of course should be avoided. To mitigate these issues, scale $\alpha_t$ to sum up to 1:
\begin{align*}
	c_t 	&	\coloneq \one^\intercal \alpha_t 	& 	\hat{\alpha}_t \coloneq \alpha_t/c_t.
\end{align*}
The update can be done in two compact steps:
\begin{equation}
\begin{split}
	\hat{\alpha}'_t &= \left(A^S(\gamma(r_{t-1}))\right)^\intercal \left( a^R_{r_{t-1}r_t}(:) \circ \hat{\alpha}_{t-1}  \right) \\
	 \hat{\alpha}_t &= \frac{\hat{\alpha}'_t}{\one^\intercal\hat{\alpha}'_t}.
\end{split}\label{eq:alpha_sc}
\end{equation}
Also note that $\log c_t = \log c_{t-1} + \log \one^\intercal\hat{\alpha}'_t$ and $\log\ell(\mu) = \log(\one^\intercal\alpha_T) = \log(c_T)$; the likelihood $\ell(\mu)$ is not computed anymore, only the log-likelihood.
Since the backwards variables $\beta_{T-t}$ can be expected to be of similar order as $\alpha_t$ they are scaled using the same scaling factors $c_t$:
\begin{equation}
	\hat{\beta}_{T-t} \coloneq \beta_{T-t}/c_{t}.
	 \label{eq:beta_sc}
\end{equation}
We then compute the scaled versions of $\xi_t(i,j)$ in similar fashion
\begin{equation}
\begin{split}
	\hat{\xi}'_t(i,j) &= \hat{\alpha}_{t-1}(i) a^R_{r_{t-1}r_t}(i) a^S_{ij}(\gamma(r_{t-1})) \hat{\beta}_t(j)\\
	\hat{\xi}_t &= \frac{ \hat{\xi}'_t}{ \one^\intercal \hat{\xi}'_t \one} = \frac{ {\xi}_t}{ \one^\intercal \alpha_T}
\end{split}\label{eq:xi_sc}
\end{equation}
and note that upon substituting $\xi_t = (\one^\intercal \alpha_T) \hat{\xi}_t$ in~\eqref{eq:reest}, $\one^\intercal\alpha_T$ cancels everywhere, and we arrive at the rescaled re-estimation equations
\begin{subequations}\label{eq:reest_sc}
	\begin{align}
		\pi^{S'}_i &= \sum_{j=1}^S \hat{\xi}_1(i,j) 	\label{eq:reestpiscale}\\ %
		a^{S'}_{ij}(l) &= \frac{\sum_{t: \gamma(r_{t-1})=l} \hat{\xi}_t(i,j)}{\sum_{k=1}^S \sum_{t: \gamma(r_{t-1})=l} \hat{\xi}_t(i,k)} 	\label{eq:reestasscale}\\
		a^{R'}_{mn}(i) &= \frac{\sum_{k=1}^S\sum_{t: r_{t-1}=m, r_t=n} \hat{\xi}_t(i,k)}{\sum_{k=1}^S \sum_{t: r_{t-1}=m} \hat{\xi}_t(i,k)}. 		\label{eq:reestarscale}
	\end{align}
	\end{subequations}
For a more detailed derivation, see again~\citep{clMMMC_arxiv}.
}
\begin{algorithm}[tbh]
\begin{algorithmic}[1]
	\Procedure{EM-Estimation}{ $\obs=(r_0\dots r_T), \mu=(\pi^R, \pi^S, A^R, A^S; \Gamma); \epsilon>0, N_{max}>0$ }
		\Statex \emph{Initialization}
			\State $N_{iter}\coloneq0$ \hfill
				 $\delta\coloneq\infty$ \hspace*{\fill} \strut
			\State $\pi^R \gets e_{r_0}$  \Comment{See remark before Lemma~\ref{thm:reest}}
		\Statex \emph{Iteration}
			\While{$N_{iter}\leq N_{max}$ and $\delta>\epsilon$}
				\Statex \emph{Forward-backward variables from $\mu$}
					\State $\begin{array}{clll}
								&\alpha \gets \eqref{eq:alpha} & \beta \gets \eqref{eq:beta} & \xi \gets \eqref{eq:xi}\\
								\text{ or }&\hat\alpha \gets \eqref{eq:alpha_sc} & \hat\beta \gets \eqref{eq:beta_sc} & \hat\xi \gets \eqref{eq:xi_sc}
						\end{array}$
				\Statex \emph{Reestimation}
					\State $(\pi^{S'}, A^{R'}, A^{S'})\gets \eqref{eq:reest} \text{ or }\eqref{eq:reest_sc}$
					\State $\mu'\gets(\pi^R, \pi^{S'}, A^{R'}, A^{S'}; \Gamma)$
				\Statex \emph{Check relative changes}
					\For{$\nu\in \{ \pi^{S}, A^{R}, A^{S}\}$}
						\State $\delta_\nu \gets \left( \|\nu-\nu'\|_F \right)/\|\nu\|_F $
					\EndFor
					\State $\delta_\ell \gets |\ell(\mu')-\ell(\mu)|/|\ell(\mu)|$
					\State $\delta \gets \max \{ \delta_\nu\}$
				\Statex \emph{Update}
					\State $N_{iter}\gets N_{iter}+1$
					\State $(\pi^{S}, A^{R}, A^{S})\gets(\pi^{S'}, A^{R'}, A^{S'})$ \hfill $\mu\gets\mu'$ \hfill\strut
			\EndWhile
			\State \textbf{return} $\mu$
	\EndProcedure
\end{algorithmic}
\caption{EM parameter estimation for cl3MCs}
\label{alg:EM}
\end{algorithm}
\section{Proofs}
\label{sec:proofs}
\begin{proof}[Lemma~\ref{lem:Q}]
Formula \eqref{eq:xi} follows directly from the definitions of $ \alpha_{t-1}$, $a^R_{mn}$, $a^S_{ij}$ and $\beta_t$.
Now let summation indices $s_\tau$ always run from 1 to $S$. First, we rewrite the second line of~\eqref{eq:Q}: given the definition of $L_{ij}  (k,\ell)$ and~\eqref{eq:Prs} we get
\begin{multline*}
    \log P(s_0,\dotsc,s_T,r_0,\dotsc,r_T\given \mu') = \log \pi^{R'}_{r_0} +     \log \pi^{S'}_{s_0} +
    \sum_{t=1}^T L_{s_{t-1} s_t}  (r_{t-1}, r_t).
\end{multline*}
We substitute this back into~\eqref{eq:Q} to get $Q(\mu,\mu') = Q_1+Q_2+Q_3$, with 
\[
Q_1= \sum_{s_0\dotsb s_T}\! P(s_0,\dotsc,s_T,r_0,\dotsc,r_T\given \mu)\cdot\allowbreak \log \pi^{R'}_{r_0} =
  \ell(\mu) \log \pi^{R'}_{r_0};
  \] for $Q_2$ and $Q_3$, note the ``marginalisation''
 \begin{multline*}
 	 \xi_t(i,j)
 		=\sum_{s_0\dots s _{t-2},s_{t+1}\dots s_T}  P\Bigl(s_0,\dotsc,s_{t-2},
	S_{t-1}=i, S_t=j,\dotsc s_T,r_0,\dotsc,r_T\given \mu\Bigr),
\end{multline*}
 where the sum runs over the states of the latent chain $\bf S$ before time instant ${t-1}$ and after time instant $t$.
 Then
 \begin{multline*}
 	 Q_2=\sum_{s_0 \dotsb s_T}P(s_0,\dotsc,s_T,r_0,\dotsc,r_T \given \mu) \log \pi^{S'}_{s_0}= \\
  		\sum_{i,j=1}^S  \log \pi^{S'}_{i} \!\sum_{s_2 \dotsb s_T}      P(S_0=i,S_1=j,s_2, \dotsc,r_T\given \mu) =
    		\sum_{i=1}^S \log \pi^{S'}_{i}\sum_{j=1}^S  \xi_1(i,j)
\end{multline*}
and
\begin{multline*}
	Q_3=
	    \sum_{s_0\dotsb s_T}
  	  	P(s_0,\dotsc,r_T\given \mu) \sum_{t=1}^T L_{s_{t-1}s_t}(r_{t-1} ,r_t)=\\
		\shoveleft{\sum_{t=1}^T\sum_{s_{t-1},s_t} L_{s_{t-1}  s_t}(r_{t-1}, r_t)\cdot}
		\shoveright{ P\bigl( S_{t-1}=s_{t-1},S_t=s_t, r_0,\dotsc,r_T \given \mu\bigr) =}\\
		 \sum_{t=1}^T\sum_{i=1}^S\sum_{j=1}^S L_{{i}j}(r_{t-1}, r_t) \xi_t(i,j).
\end{multline*}
This completes the proof.
\end{proof}
\begin{proof}[Theorem~\ref{thm:reest}]
	From the remark before the theorem, it should be clear that we
        can ignore the first term in~\eqref{eq:Qsum} in the
        maximisation. Consider $\mu\in\Lambda$ as fixed, and define
        $\widetilde\mu\mapsto W(\widetilde\mu):=Q(\mu,\widetilde\mu)$.
        We claim that $W$ has the unique global maximum point
        $\mu'\in\Lambda$.
         Note that if $\mu\in\Lambda$ then it may have zero components, say $\mu_{i_1}=\dots=\mu_{i_d}=0$. Then the logarithms of the corresponding components $\tilde{\mu}_{i_1}\dots\tilde{\mu}_{i_d}$ of $\tilde\mu$ in $W$ are multiplied by $0$ so that these components do not change $W$. However, if we fix all the components of $\tilde\mu$ but $\tilde\mu_{i_1}$ and any other component $\tilde\mu_k$ such that $k\not\in\{i_1\dots i_d\}$, and $k$ is such that $\tilde\mu_{i_1}$ and $\tilde\mu_k$ are in $\pi^R$, or are in the same row of $A^R$ or $A^S$, then increasing $\tilde\mu_{i_1}$ will decrease $\tilde\mu_k$ (to meet the stochasticity constraints). As a result, the $\log \tilde\mu_k$ will decrease causing $W$ to decrease.
         Hence, the maximum of $W$ is attained in the set $\widetilde\Lambda\coloneq\{\tilde\mu\in\Lambda: \tilde{\mu}_{i_1}=\dots=\tilde{\mu}_{i_d}=0\}$. Let
$\widetilde W$ be the restriction of $W$ to $\widetilde\Lambda$. Now $\widetilde W$ is a conical sum of logarithms of \emph{all} $R+S+SR^2+pS^2-d$ independent components of $\tilde\mu$, hence \emph{strictly} concave function on a convex compact set $\widetilde\Lambda$. Hence, $\widetilde W$ has a unique maximum point in $\widetilde \Lambda$ which coincides with the unique global maximum point $\mu'$ of $W$ in $\Lambda$.

Let us prove that $\mu'=M(\mu)$. As noted above, $\mu_{i_1}=\dots=\mu_{i_d}=0$ implies that $\mu'_{i_1}=\dots=\mu'_{i_d}=0$, and we stress that the same property holds true for $M(\mu)$: %
  If $\pi^S_i$ is $0$ it follows from~\eqref{eq:alpha} that $\alpha_0(i)=0$ and from~\eqref{eq:xi} we get $\xi_1(i,j)=0$. By~\eqref{eq:reestpi}, $\pi^{S'}_i=0$ as well. If we have $a^S_{ij}(l)=0$, then again from~\eqref{eq:xi}, we see that $\xi_t(i,j)=0$ for all $t$ with $r_{t-1}=l$, and~\eqref{eq:reestas} yields $a_{ij}^{S'}(l)=0$. Similarly, $a^R_{mn}(i)=0$ leads to $\xi_t(i,j)=0$ whenever $r_{t-1}=m$ and $r_t=n$, independently of $j$, so $a^{R'}_{mn}(i)=0$, too. Hence $M(\widetilde\Lambda)=\widetilde\Lambda$, and $\pi^{S'}_i$, $a^{S'}_{ij}(l)$ and $a^{R'}_{mn}(i)$ defined
by~\eqref{eq:reest} are positive if the corresponding components of $\mu$ are. On the other hand, $\mu'_k>0$ if $\mu_k>0$ as otherwise $W(\mu')=-\infty$, and so the gradient of $\widetilde W$ is well-defined at $\mu'$. In fact, for positive $\tilde{\pi}^S_i$, $\tilde{a}^S_{ij}(l)$ and $\tilde{a}^R_{mn}(i)$:
\begin{align*}
	\ddpart{\widetilde W}{\tilde{\pi}^S_i} &= \frac{\sum_{j=1}^S\xi_1(i,j)}{\tilde{\pi}^S_i}\\
	\ddpart{\widetilde W}{\tilde{a}^S_{ij}(l)} & = \frac{\sum_{t: \gamma(r_{t-1})=l} \xi_t(i,j)}{\tilde{a}^S_{ij}(l)}\\
	\ddpart{\widetilde W}{\tilde{a}^R_{mn}(i)} & = \frac{\sum_{k=1}^S\sum_{t:r_{t-1}=m, r_t=n}\xi_t(i,k)}{\tilde{a}^R_{mn}(i)}
\end{align*}
By e.g.~\citep[p.113, Prop. 2.1.2]{Bertsekas1999}, it is necessary and sufficient for $\mu'$ to satisfy the inequality
\begin{equation}
  \label{eq:support}
(\nabla \widetilde W(\widetilde\mu)\Bigr|_{\widetilde\mu=\mu'})^\intercal\widetilde\mu\le
(\nabla \widetilde W(\widetilde\mu)\Bigr|_{\widetilde\mu=\mu'})^\intercal\mu' \quad \forall\widetilde\mu\in\widetilde\Lambda .
\end{equation}
We stress that the r.h.s.\ of~\eqref{eq:support} is independent of $\mu'$ as
\begin{equation*}
    (\nabla \widetilde W(\widetilde\mu)\Bigr|_{\widetilde\mu=\mu'})^\intercal\mu' =(2T+1) \ell(\mu).
\end{equation*}
Since $\pi^{S'}_i$, $a^{S'}_{ij}(l)$ and $a^{R'}_{mn}(i)$ defined
by~\eqref{eq:reest} are positive if the corresponding components of
$\mu$ are so, the gradient of $\widetilde W$ is well-defined at
$M(\mu)$. Take any $\widetilde\mu\in\widetilde\Lambda$ and compute:
\begin{equation*}
  \begin{split}
    \widetilde{\pi}^{S}_i \ddpart{\widetilde W}{\tilde{\pi}^S_i}(\mu')
    &= \widetilde{\pi}^{S}_i\sum_{i=1}^S\sum_{j=1}^S\xi_1(i,j)\\
    \widetilde{a}^{S}_{ij}(l) \ddpart{\widetilde W}{\tilde{a}^S_{ij}(l)}(\mu') &= \widetilde{a}^{S}_{ij}(l)\sum_{k=1}^S \sum_{t:
      \gamma(r_{t-1})=l} \xi_t(i,k)\\
    \widetilde{a}^{R}_{mn}(i) \ddpart{\widetilde W}{\tilde{a}^R_{mn}(i)}(\mu')&=\widetilde{a}^{R}_{mn}(i)\sum_{n=1}^R\sum_{k=1}^S
                                 \sum_{t: r_{t-1}=m, r_t=n}                            \xi_t(i,k)
  \end{split}
\end{equation*}
so that, by stochasticity constraint, we get~\citep{clMMMC_arxiv}:
\begin{equation*}
    (\nabla \widetilde
W(\widetilde\mu)\Bigr|_{\widetilde\mu=\mu'})^\intercal\widetilde{\mu} = (2T+1) \ell(\mu)
\end{equation*}
Hence, $M(\mu)$ defined by~\eqref{eq:reest}
satisfies~\eqref{eq:support} with equality for any
$\widetilde{\mu}\in\widetilde{\Lambda}$. This completes the proof.
\end{proof}

\section{Experimental details}\label{app:exp}
The detailed steps in generating the data in Section~\ref{ssec:synth} are as follows.
We consider ${T}\in\{500 ,\allowbreak1000,\allowbreak3000,\allowbreak5000 ,\allowbreak7500 ,\allowbreak10000 ,\allowbreak12500 ,\allowbreak15000 ,\allowbreak 20000 ,\allowbreak25000 ,\allowbreak30000 ,\allowbreak50000 ,75000 \}$ and begin at $i=0$.  Then, for each of the $N_T=13$ values of $T$, the following is repeated $N_e$ times:
\begin{enumerate*}[start=1]
	\item generate a pair of $R\times R$ row-stochastic matrices $\left(A^R(1), A^R(2)\right)$ by selecting one ``dominant element'' per row and setting it to a random number uniformly distributed in $[0.5,1]$ (this is done to ensure that the two matrices in a pair are sufficiently different) and filling the remaining elements with uniformly random numbers;
	\item generate a $S\times S$ row-stochastic matrix $A^S(1)$ with uniformly random entries;
	\item generate initial probability vectors $\pi^S\in [0,1]^S$ and $\pi^R \in [0,1]^R$ with uniformly random entries;
	\item generate a trajectory of length $T$ from the cl3MC $\mu_o = (\pi^S, \pi^R, A^S(1), A^R; [R] )$ (this is the open-loop case);
	\item estimate parameters of $\mu_o$ by running Algorithm~\ref{alg:EM} (initialized with uniformly random guesses for unknown parameters) until convergence and record the distance measures defined in Section~\ref{ssec:synth};
	\item generate an additional $S\times S$ row-stochastic matrix $A^S(2)$ with uniformly random entries;
	\item generate a random partition $\Gamma$ of $[R]$ by first randomly permuting $[R]$ and then splitting it after a random index between $1$ and $R-2$;
	\item generate a trajectory of length $T$ from the cl3MC $\mu_c = (\pi^S, \pi^R, A^S, A^R; \Gamma )$ (this is the open-loop case);
	\item estimate parameters of $\mu_c$ by running Algorithm~\ref{alg:EM} (initialized with uniformly random guesses for unknown parameters) until convergence and record the distance measures defined in Section~\ref{ssec:synth}.
\end{enumerate*}
Note that we assume $A^R(2)$ to be known and it is not estimated. Concretely that means that the initial guess of $A^R(2)$ is set to the true value of $A^R(2)$ and $A^R(2)$ is excluded from the reestimation steps in Algorithm~\ref{alg:EM}. The parameters used in Algorithm~\ref{alg:EM} are $\varepsilon=10^{-5}$, a typical, empirical choice for relative tolerances, and $N_{max}=2000$, at which point the algorithm has typically long converged (in our experiments, $N_{max}$ was never reached).
 Once the experiment terminates, we have $2 N_e N_T=2600$ (resp.\  $3 N_e N_T=3900$) different values for each of the distance metrics for $A^R$ (resp.\ $A^S$), which are visualized in Figures~\ref{fig:err_vs_T} and \ref{fig:page_fixed_bars_AR}.

For Section~\ref{ssec:driver}, the matrices $\left(A^R(1), A^R(2)\right)$ are fixed to the ones indicated in Figure~\ref{fig:schoolmap}. One run then consisted of repeating $N_t$ times the following:
\begin{enumerate*}
\item draw a number from the uniform distribution on $[0,1]$ and if it is less than $p$, set $s=1$, else $s=2$;
\item using the transition probabilities in $A^R(s)$, generate a trajectory from initial state $O$ until state $D$ is reached.
\end{enumerate*}
We then concatenated all of those trajectories, identifying states $O$ and $D$ as state 1, e.g.\ if there were only two trajectories, $(O,2,6,5,9,D)$ and $(O,7,4,5,8,D)$ they are concatenated to $(1,2,6,5,\allowbreak9,1,7,4,5,8,1)$. Algorithm~\ref{alg:EM} is then run on that single trajectory with the same parameters as above, again assuming $A^R(2)$ as known, two times: once assuming open-loop, i.e.\ $\Gamma=[9]$, and once assuming closed-loop with $\Gamma = \bigl\{ \{8,9\}, \{1,\dotsc,7\} \bigr\}$ and $A^S(2)=I_2$. Then the distance metrics between estimates and true values are computed. This concluded one run; we collected $N_e$ runs, yielding $2 N_e=100$ numbers for each $A^R(1)$ and $A^S(1)$. The numbers for $A^R(1)$ are shown in Figure~\ref{fig:driver_results}.

\end{document}